\documentclass[12pt]{amsart}%
\usepackage{amsfonts}
\usepackage{amsmath}
\usepackage{amssymb}
\usepackage{graphicx}%
\newtheorem{theorem}{Théorème}
\theoremstyle{plain}

\newtheorem{corollary}{Corollaire}

\newtheorem{definition}{Définition}
\newtheorem{example}{Exemple}

\newtheorem{lemma}{Lemme}
\newtheorem{notation}{Notation}

\newtheorem{proposition}{Proposition}
\newtheorem{remark}{Remarque}

\numberwithin{equation}{section}
\begin{document}
\title{Th\'{e}orie De Galois Des Equations Diff\'{e}renti\`{e}lles}
\author{Lotfi Sa\"{\i}dane}
\address{LOTFI\ SA\"{I}DANE\\
D\'{e}partement de Math\'{e}matiques \\
Facult\'{e} des sciences de Tunis\\
Universit\'{e} De Tunis El-Manar\\
2092 El-ManarI, TUNIS}
\email{Lotfi.saidane@fst.rnu.tn}
\date{15 Decembre 2002}

\begin{abstract}
Let $k$ be a differential field and $C$ its subfield of constants. In general
a differential extension $K$ of $k$ add some new constants to $C,$ and it is
difficult to prove that $C$ stay unchangeable under the extension $K;$ This
situation is provided by the Picard-Vessiot extension. Kolchin \cite{Kolchin}
prove the theorem of existence and unicity for these extensions. The aim of
this paper is to prove Kolchin theorem and other results, in a simple manner,
by means of the theory of models and logic.

\bigskip\medskip

\noindent\textbf{R\'{E}SUM\'{E}: }Soit $k$ un corps diff\'{e}rentiel et $C$
son sous corps des constantes. En g\'{e}n\'{e}ral une extension
diff\'{e}rentiel $K$ de $k$ modifie le corps des constantes $C$ de $k.$
Prouver que $K$ ne modifie pas $C$ est un probl\`{e}me assez difficile en
alg\`{e}bre diff\'{e}rentiel. Les extensions de Picard-Vessiot constitue un
exemple de cette situation. Kolchin \cite{Kolchin} a montr\'{e} le
th\'{e}or\`{e}me d'existence et d'unicit\'{e}, \`{a} isomorphisme pr\'{e}s,
des extensions de Picard-Vessiot sous la condition que le corps $C$ est
alg\'{e}briquement clos. Dans ce travail on utilise la th\'{e}orie des corps
diff\'{e}rentiellements clos (Th\'{e}orie des mod\`{e}les), pour montrer
l'existence et l'unicit\'{e}, \`{a} isomorphisme pr\'{e}s, des extensions de
Picard-Vessiot. On calcul ensuite le groupe de Galois diff\'{e}rentiel de
certaines extensions particuli\`{e}re. Enfin, on montre quelque
th\'{e}or\`{e}mes g\'{e}n\'{e}raux de la th\'{e}orie de Galois
diff\'{e}rentielle par les m\^{e}mes techniques.

\end{abstract}
\maketitle

\section{Notations\textbf{{\protect\Large \ }}Et\textbf{{\protect\Large \ }%
}Pr\'{e}liminaires}

Un anneau diff\'{e}rentiel, est un anneau commutatif $A$, muni d'une
d\'{e}rivation, c'est \`{a} dire d'une application $\delta:A\rightarrow A$
satisfaisant pour tout $x$ et $y$ dans $A,$ $\delta(x+y)=\delta x+\delta y$ et
$\delta(xy)$= $x\delta y+y\delta x.$ Par exemple, l'anneau des fonctions
m\'{e}romorphes sur $\mathbb{C}$ muni de la d\'{e}rivation $\partial/\partial
z;$ l'anneau des s\'{e}ries formelles sur un corps $k$ (quelconque), $k[[x]]$
muni de la d\'{e}rivation $\partial/\partial x;$ sont des anneaux
diff\'{e}rentiels. On appelle constante de $A$, les \'{e}l\'{e}ments de $A$ de
d\'{e}riv\'{e}e nulle. Les constantes forment un sous anneau de $A.$ On
appelle id\'{e}al diff\'{e}rentiel de $A$, tout id\'{e}al $I$ de $A$
v\'{e}rifiant $\delta I$ $\subset I.$ Si $A$ est anneau diff\'{e}rentiel
int\`{e}gre sa d\'{e}rivation s'\'{e}tend de mani\`{e}re unique \`{a} son
corps des fractions. Naturellement un corps diff\'{e}rentiel est un corps muni
d'une d\'{e}rivation. Les constantes d'un corps diff\'{e}rentiel constituent
un corps. Pour la suite de ce travail, on ne consid\`{e}re que des corps
diff\'{e}rentiels de caract\'{e}ristique z\'{e}ro. Soit $k$ un corps
diff\'{e}rentiel, on note par $k\{x\}$ l'anneau $k[x,$ $x^{\prime},$
$x^{\prime\prime},$ $...,$ $x^{(n)},$ $...]$, o\`{u} $n\in\mathbb{N},$ des
polyn\^{o}mes en une infinit\'{e} d\'{e}nombrable de variables, muni de
l'unique d\'{e}rivation prolongeant celle de $k$, qu'on note $\partial,$
v\'{e}rifiant pour tout $n\in\mathbb{N}$, $\partial x^{(n)}=x^{(n+1)}.$ On
appelle "type" en une variable \`{a} coefficient dans $k$, tout id\'{e}al
diff\'{e}rentiel premier de $k\{x\}$ distinct de $k\{x\},$ ce dernier est
exclu de l'ensemble des id\'{e}aux diff\'{e}rentiels premiers. Soient
$K\supset k$ une extension de corps diff\'{e}rentiel, $a$ et $b$ deux
\'{e}l\'{e}ment de $K.$ On dit que $a$ et $b$ ont "m\^{e}me type" sur $k,$
s'il existe un isomorphisme entre $k\{a\}$ et $k\{b\}$ qui fixe les
\'{e}l\'{e}ments de $k$ et associe $\partial^{n}a$ \`{a} $\partial^{n}b$
($n\in\mathbb{N}$). L'application $\varphi:k\{x\}\rightarrow k\{a\},$ qui
\`{a} $P(x)$ associe $P(a)$ est un homomorphisme d'anneau diff\'{e}rentiel,
son noyau est un id\'{e}al diff\'{e}rentiel premier, qu'on note par $I_{a}$
(id\'{e}al des \'{e}quations \`{a} coefficients dans $k$ satisfaites par $a)$.
On a, $k\{a\}\simeq k\{x\}/I_{a}$ et le corps des fractions de $k\{a\},$ soit
$k\langle a\rangle,$ est isomorphe au corps des fractions de $k\{x\}/I_{a}.$
Deux \'{e}l\'{e}ments $a$ et $b$ dans $K\supset k$, ont alors "m\^{e}me type"
sur $k$ si et seulement si $I_{a}=I_{b}.$ L'id\'{e}al $I_{a}$ est appel\'{e}
type de $a.$ Tout id\'{e}al diff\'{e}rentiel premier de $k\{x\}$ est le type
d'un \'{e}l\'{e}ment $a$ dans une extension diff\'{e}rentielle de $k.$
L'ensemble des types en une variable \`{a} coefficients dans $k$ ( ou encore
$1-type$) est not\'{e} $S_{1}(k).$

\subsection{ Topologie de $S_{1}(k)$}

Soit $P$ un \'{e}l\'{e}ment de $k\{x\},$ on note par $\langle P\rangle$
l'id\'{e}al, au sens habituel, engendr\'{e} par $P$ et ses
d\'{e}riv\'{e}\'{e}s successives, c'est le plus petit id\'{e}al
diff\'{e}rentiel contenant $P$. Plus g\'{e}n\'{e}ralement, si $\mathcal{P}$
est un ensemble de polyn\^{o}mes diff\'{e}rentiels de $k\{x\},$ l'id\'{e}al
$\langle\mathcal{P\rangle}$ est le plus petit id\'{e}al diff\'{e}rentiel
contenant $\mathcal{P}$.

\begin{definition}
On appelle ordre de $P\in k\{x\}\backslash k,$ le plus grand ordre de la
d\'{e}riv\'{e}e de $x$ qui intervient dans $P.$
\end{definition}

Par convention, le polyn\^{o}me nul n'a pas d'ordre et les \'{e}l\'{e}ments de
$k^{\ast}$sont d'ordre $-1.$

\begin{definition}
Soit $P\in k\{x\}$ d'ordre $n\geq0,$ on appelle s\'{e}parante de $P$, et on
note $S_{P}$ le polyn\^{o}me $\delta P/\delta x^{(n)}.$
\end{definition}

La s\'{e}parante de $P$ est un polyn\^{o}me non nul et d'ordre inf\'{e}rieur
o\`{u} \'{e}gal \`{a} celui de $P.$ Si $P$ et $S_{P}$ sont d'ordre \'{e}gal
\`{a} $N$, alors le degr\'{e} partiel en $x^{(N)}$ dans $S_{P}$ est
strictement inf\'{e}rieur \`{a} celui de $P.$ D'o\`{u}, un polyn\^{o}me ne
divise pas sa s\'{e}parante. Soit $P$ un polyn\^{o}me diff\'{e}rentiel d'ordre
$N\geq0,$ sans facteur commun avec sa s\'{e}parante, alors l'id\'{e}al
diff\'{e}rentiel $\langle P\rangle$ ne contient pas de polyn\^{o}mes
diff\'{e}rentiels d'ordre strictement inf\'{e}rieur \`{a} $N$ et ses
\'{e}l\'{e}ments d'ordre $N$ sont multiples de $P,$ voir \cite{Poizat78,
Poizat77} pour une d\'{e}monstration de ce r\'{e}sultat et des deux Lemmes qui
vont suivrent.

\begin{lemma}
Soit $P$ un polyn\^{o}me diff\'{e}rentiel irr\'{e}ductible d'ordre $N\geq0.$
On consid\`{e}re l'ensemble $I(P)$, des plyn\^{o}mes diff\'{e}rentiels $Q$
v\'{e}rifiant pour $m$ assez grand ($S_{P})^{m}Q\in\langle P\rangle.$ Alors
$I(P)$ est un id\'{e}al diff\'{e}rentiel premier contenant $P$ et ne contenant
aucun \'{e}l\'{e}ment d'ordre strictement inf\'{e}rieur \`{a} $N.$ Les
\'{e}lements d'ordre $N$ de $I(P)$ sont multiple de $P.$

\begin{proof}
[Preuve]Pour la d\'{e}monstration voir \cite{Poizat78}
\end{proof}
\end{lemma}

\begin{lemma}
Tout id\'{e}al diff\'{e}rentiel premier non nul est de la forme $I(P).$

\begin{proof}
[Preuve]Pour la d\'{e}monstration voir \cite{Poizat78}
\end{proof}
\end{lemma}

\begin{proposition}
L'id\'{e}al $I(P)$ est le plus petit id\'{e}al diff\'{e}rentiel premier
contenant $P$ et ne contenant pas $S_{P}.$

\begin{proof}
[Preuve]Soit $I$ un id\'{e}al diff\'{e}rentiel premier contenant $P$ et ne
contenant pas $S_{P}.$ Soit $Q$ $\in I(P)$, il existe alors $m$ dans
$\mathbb{N}$ tel que $S_{P}^{m}Q\in\langle P\rangle,$ comme $S_{P}\notin I$
$\supset\langle P\rangle,$ alors $Q\in I$ et $I(P)\subset I.$
\end{proof}
\end{proposition}

\begin{example}
Si $P=(X^{\prime})^{2}-2X$ est dans $k\{x\},$ $P$ est irr\'{e}ductible
ind\'{e}pendamment du corps $k$ choisi. Sa s\'{e}parante $S_{P}=2X^{\prime}$.
En vertue de la discussion pr\'{e}c\'{e}dant le Lemme1, l'id\'{e}al
diff\'{e}rentiel $\langle P\rangle$ ne contient pas $S_{P}.$ Tout id\'{e}al
diff\'{e}rentiel qui contient $P$ et $S_{P}$ contient $X,$ le plus petit
\'{e}tant $\langle X\rangle$ ($I_{X}=\langle X\rangle).$ On remarque que
l'id\'{e}al diff\'{e}rentiel $\langle X\rangle\nsubseteq I_{P},$ car sinon
$X^{\prime}$serait dans $I_{P}.$ Le polyn\^{o}me d\'{e}riv\'{e}e de $P,$
$P^{\prime}=2X^{\prime}X"-2X^{\prime}=$ $2X^{\prime}(X"-1),$ comme $I_{P}$ est
premier et $X^{\prime}\notin I_{P}$ alors $X"-1$ $\in I_{P}$ mais
$X"-1\notin\langle X\rangle$ et $I_{P}\nsubseteq\langle X\rangle.$
L'intersection des id\'{e}aux diff\'{e}rentiels $I_{P}$ et $\langle X\rangle,$
\'{e}tant un id\'{e}al diff\'{e}rentiel contenant $\langle P\rangle$ et ne
contenant pas $S_{P},$ de la proposition 1 on d\'{e}duit que $\langle
P\rangle=I_{P}\cap\langle X\rangle.$
\end{example}

En g\'{e}n\'{e}ral, avec la condition d'irr\'{e}ductibilit\'{e}e, les
id\'{e}aux diff\'{e}rentiels $I_{Q}$ et $\langle Q\rangle,$ sont identiques.
L'exemple pr\'{e}c\'{e}dent \'{e}tablit une situation de diff\'{e}rence.

Les id\'{e}aux diff\'{e}rentiels premiers sont de la forme $I_{P},$ o\`{u} $P$
est un polyn\^{o}me diff\'{e}rentiel irr\'{e}ductible (Lemme 2). On a une
bijection entre $S_{1}(k)$ et l'ensemble des polyn\^{o}mes diff\'{e}rentiels irr\'{e}ductibles.

\begin{proposition}
Tout id\'{e}al diff\'{e}rentiel premier qui contient un polyn\^{o}me
diff\'{e}rentiel irr\'{e}ductible $P$ d'ordre $N\geq0$ et qui n'est pas
$I_{P}$ contient un \'{e}l\'{e}ment d'ordre strictement inf\'{e}rieur \`{a}
$N.$
\end{proposition}

\begin{proof}
[Preuve]Soit $I$ un id\'{e}al diff\'{e}rentiel premier, d'apr\`{e}s le lemme
2, il existe alors un polyn\^{o}me diff\'{e}rentiel irr\'{e}ductible $Q$ tel
que $I=I_{Q}$. Comme l'ordre de $Q$ est inf\'{e}rieur o\`{u} \'{e}gal \`{a}
$N,$ car $P\in I$ et le lemme 1. On Suppose que l'ordre de $Q$ est \'{e}gal
\`{a} $N$, le polyn\^{o}me $Q$ diviserait alors $P,$ ce qui contredit le fait
que $P$ est irr\'{e}ductible.
\end{proof}

\begin{remark}
Soit $a$ un \'{e}l\'{e}ment d'une extension diff\'{e}rentielle $K$ de $k.$
Supposons que le type de $a$ sur $k$ est $I_{P},$ o\`{u} $P$ est un
polyn\^{o}me diff\'{e}rentiel irr\'{e}ductible d'ordre $N\geq0.$ Alors $N$ est
le degr\'{e} de transcendance au sens alg\'{e}brique du corps $k(a)$ sur $k.$
\end{remark}

Soit $P$ un polyn\^{o}me diff\'{e}rentiel irr\'{e}ductible de $k\{x\},$ on
consid\`{e}re l'ensemble $\mathcal{O}_{P}=\{I\in S_{1}(k)$ $/$ $P\in I\}.$ On
munit $S_{1}(k)$ de la topologie, appell\'{e}e topologie de Stone, dont une
base d'ouverts est constitu\'{e}e par l'ensemble des $\mathcal{O}_{P}$.

\begin{proposition}
L'ensemble $S_{1}(k)$ est compact et les ouverts ferm\'{e}s de $S_{1}(k)$ sont
exactement les $\mathcal{O}_{P}.$

\begin{proof}
[Preuve]Voir \cite{Ehrsam}
\end{proof}
\end{proposition}

\begin{definition}
soit $K$ un corps diff\'{e}rentiel, $K$ est dit diff\'{e}rentiellement clos si
tout syst\`{e}me de la forme:\qquad$\{%
\begin{array}
[c]{c}%
P(X)=0\\
Q(X)\neq0
\end{array}
$\newline o\`{u} l'ordre de $Q$ est strictement inf\'{e}rieur \`{a} celui de
$P,$ \`{a} coefficients dans $K$ a une solution dans $K.$
\end{definition}

\begin{example}
$\{%
\begin{array}
[c]{cc}%
X"X-X^{\prime} & =0\\
X^{\prime} & \neq0
\end{array}
$
\end{example}

On remarque qu'un corps diff\'{e}rentiellement clos est ag\'{e}briquement
clos. Il suffit de consid\'{e}rer un syst\`{e}me constitu\'{e} par une
\'{e}quation polyn\^{o}miale, c'est \`{a} dire d'ordre $0$ et une
in\'{e}quation constante, d'ordre $-1.$

\begin{lemma}
Tout corps diff\'{e}rentiel se plonge dans un corps diff\'{e}rentiellement clos.
\end{lemma}

\begin{definition}
Le wronskien de $n$ \'{e}l\'{e}ments $y_{1},...,y_{n}$ dans un anneau
diff\'{e}rentiel est le d\'{e}terminant de la matrice de wronsky
\[
\left\vert
\begin{array}
[c]{llll}%
y_{1} & y_{2} & \cdots & y_{n}\\
y_{1}^{\prime} & y_{2}^{\prime} & \cdots & y_{n}^{\prime}\\
y_{1}^{(n-1)} & y_{2}^{(n-1)} & \cdots & y_{n}^{(n-1)_{n}}%
\end{array}
\right\vert
\]

\end{definition}

\begin{theorem}
Soit $K$ un corps diff\'{e}rentiel, son corps des constantes \'{e}tant $C$;
$n$ \'{e}l\'{e}ments de $K$ sont lin\'{e}airement d\'{e}pendants sur $C$ si et
seulement si leurs wronskien est nul.

\begin{proof}
Supposons que $y_{1},...,y_{n}$ sont lin\'{e}airement d\'{e}pendants sur $C.$
Soit $C_{1},$ ..., $C_{n},$ non tous nul, dans $C$ tel que $\sum_{i=1}%
^{n}C_{i}y_{i}=0$. En d\'{e}rivant cette \'{e}quation $(n-1)$ fois, on obtient
alors $n$ \'{e}quations en $(C_{1},...,C_{n})$ de la forme $\sum_{i=1}%
^{n}C_{i}y_{i}^{(j)}=0,\quad j=0,1,...,$ $n-1.$ Comme ce syst\`{e}me admet une
solution non triviale $(C_{1},...,C_{n})$ on d\'{e}duit que le wronskien de
$(y_{1},y_{2},...,y_{n})$ est nul. Pour prouver la r\'{e}ciproque, supposons
que le wronskien de $(y_{1},$ $y_{2},$ ..., $y_{n})$ est nul. On peut alors
trouver dans $K^{n}$ une solution non triviale $(C_{1},C_{2},...,C_{n})$ du
syst\`{e}me $\sum_{i=1}^{n}C_{i}y_{i}^{(j)}=0,\quad j=0,...,n-1.$ Sans pertes
de g\'{e}n\'{e}ralit\'{e}s, on peut supposer que $C_{1}=1$ et le wronskien de
$(y_{2},...,y_{n})$ non nul. En effet, $K$ est un corps, on peut multiplier
par l'inverse de $C_{1}$ (suppos\'{e} non nul) le syst\`{e}me d'\'{e}quations.
Si le wronskien de $(y_{2},...,y_{n})$ est nul, on peut trouver dans $K^{n-1}$
une solution non triviale $(C_{2},...,C_{n})$ du syst\`{e}me $\sum_{i=2}%
^{n}C_{i}y_{i}^{(j)}=0,$ $j=0,$ $1,$ ...., $n-2.$ On peut supposer que
$C_{2}=1,$ pour les m\^{e}mes raisons que pr\'{e}c\'{e}demment, le wronskien
de $(y_{3},...,y_{n})$ est soit nul, soit non nul. Si pour tout $k=1,...,n-1$
le wronskien de $(y_{k},...,y_{n})$ est nul$,$ par le proc\'{e}d\'{e}
d\'{e}j\`{a} utilis\'{e}, on obtient des solutions $\{C_{i}\}$ $i=1,$ $2,$
$...,$ $n-1$ du syst\`{e}me qui sont \'{e}gales soit \`{a} z\'{e}ro soit \`{a}
1, ainsi que le syst\`{e}me:
\[
\left\{
\begin{array}
[c]{ll}%
y_{n-1}+c_{n}y_{n} & =0\\
y_{n-1}^{\prime}+c_{n}y_{n}^{\prime} & =0
\end{array}
\right.
\]
avec le wronskien de $(y_{n-1},y_{n})$ nul donc $y_{n-1}y_{n}^{\prime}%
-y_{n-1}^{\prime}y_{n}=0.$ Comme $C_{n}=-\frac{y_{n-1}}{y_{n}}$ on a
$C_{n}^{\prime}=-\frac{y_{n-1}^{\prime}y_{n}-y_{n-1}y_{n}^{\prime}}{y_{n}^{2}%
}=0$ donc $C_{n}$ est une constante, et le syst\`{e}me $(y_{1},...,y_{n})$ est
lin\'{e}airement d\'{e}pendant sur les constantes.\newline S'il existe
$k=1,...,n-1$ tels que le wronskien de $(y_{k},...,y_{n})$ soit non nul, on
d\'{e}rive $(n-1)$ fois l'\'{e}quation $y_{1}+C_{2}y_{2}$ $+...+$
$C_{k-1}y_{k-1}+...+$ $C_{n}y_{n}$ o\`{u} $C_{2},$ $...,$ $C_{k-1}$ sont soit
z\'{e}ro soit 1. De l'hypoth\`{e}se sur le wronskien de $(y_{1},$ $y_{2},$
..., $y_{n}),$ on d\'{e}duit l'existence d'une solution non triviale $(C_{2},$
...., $C_{n})$ du syst\`{e}me homog\`{e}ne $y_{1}^{(j)}+\sum_{i=2}^{n}%
C_{i}y_{i}^{(j)}=0,$ $j=0,$ $...,$ $n-1.$ D'o\`{u}, apr\'{e}s simplification
du syst\`{e}me initial on obtient, $\sum_{i=k}^{n}C_{i}^{\prime}y_{i}^{(j)}%
=0$, $j=0,...,$ $k-1$ comme le wronskien de $(y_{k},$ $...,$ $y_{n})$ est non
nul, on deduit que $C_{k}^{\prime}=$ $C_{k+1}^{\prime}=...$ $=$ $C_{n}%
^{\prime}=0$ et les $C_{i}$ sont des constantes.
\end{proof}
\end{theorem}

\section{Extension de Picard-Vessiot}

On consid\`{e}re l'\'{e}quation diff\'{e}rentielle lin\'{e}aire homog\`{e}ne
\[
(\ast)\qquad y^{(n)}+a_{1}y^{(n-1)}+...+a_{n-1}y^{\prime}+a_{n}y=0
\]
o\`{u} les coefficients $a_{i}$ sont dans un corps diff\'{e}rentiel
$K$.\newline Soit $u_{1},...,u_{n}$ des solutions de l'\'{e}quation dans une
certaine extension de $K$. On dit que $(u_{1},...,u_{n})$ forme un syt\`{e}me
fondamental de solutions si leurs wronskien est non nul.

\begin{definition}
Soit $y^{(n)}+a_{1}y^{(n-1)}+...+a_{n-1}y^{\prime}+a_{n}y=0$ une \'{e}quation
diff\'{e}rentielle lin\'{e}aire, les coefficients $a_{i}$ sont dans
$K$.\newline On dit que le corps $M\supseteq K$ est une extension de
Picard-Vessiot de $K$ si :\newline1) $M=K(u_{1},...,u_{n})$ o\`{u}
$(u_{1},...,u_{n})$ est un syst\`{e}me fondamental de solution.\newline2) $M$
a le m\^{e}me corps de constantes que $K$
\end{definition}

\subsection{Existence des extensions de Picard-Vessiot}

Pour montrer l'existence des extensions de Picard-Vessiot on plonge $K$ dans
un corps diff\'{e}rentiellement clos $L$ (on peut prendre $L=M(k)$ cl\^{o}ture
diff\'{e}rentielle de $K,$ voir lemme 3 et pour plus de d\'{e}tails
\cite{Poizat78} chapitre IV). On consid\`{e}re le syst\`{e}me:
\[
\left\{
\begin{array}
[c]{ll}%
y^{(n)}+a_{1}y^{(n-1)}+...+a_{n}y & =0\\
y^{\prime} & \neq0
\end{array}
\right.
\]
Les coefficients sont dans $L$, on suppose que $n>1$ (le cas $n=1$ est
trivial); puisque $L$ est un corps diff\'{e}rentiellement clos. Ce syst\`{e}me
admet une solution dans $L$ not\'{e}e $u_{1}$. \newline On consid\`{e}re
ensuite le syst\`{e}me
\[
\left\{
\begin{array}
[c]{ll}%
y^{(n)}+a_{1}y^{(n-1)}+...+a_{n}y & =0\\
y^{\prime}u_{1}-u_{1}^{\prime}y & \neq0
\end{array}
\right.
\]
Pour les m\^{e}mes raisons, ce syst\`{e}me admet une solution $u_{2}$ dans
$L$. Par induction on obtient les solutions $u_{1},u_{2},...,u_{m-1},\quad
m\leq n.$ On consid\`{e}re, alors, le syst\`{e}me
\[
\left\{
\begin{array}
[c]{ll}%
y^{(n)}+a_{1}y^{(n-1)}+...+a_{n-1}y^{\prime}+a_{n}y & =0\\
w(y,u_{1},...,u_{m-1}) & \neq0
\end{array}
\right.
\]
Ce syst\`{e}me admet une solution not\'{e}e $u_{m},$ puisque $w(y,u_{1}%
,...,u_{m-1})$ est une \'{e}quation diff\'{e}rentielle lin\'{e}aire d'ordre
$(m-1)$ en $y$. Comme $m-1<n$ et $L$ est diff\'{e}rentiellement clos, le
syst\`{e}me pr\'{e}c\'{e}dent admet une solution not\'{e}e $u_{m}$ dans $L$.
$(u_{1},...,u_{n})$ est un syst\`{e}me fondamental de solutions puisque
$w(u_{1},...,u_{n})\neq0$. Montrons que dans la cl\^{o}ture
diff\'{e}rentielle, l'extension $K(u_{1},...,u_{n})$ est de Picard-Vessiot.
Pour cela regardons le corps des constantes de $K(u_{1},...,u_{n})$.\newline
On consid\`{e}re l'ouvert ferm\'{e} pour la topologie de Stone de la forme
$\langle x^{\prime}=0\rangle$ (voir Proposition 3 et pour plus de d\'{e}tails
\cite{Ehrsam} Proposition 1, page1-04), il contient d'une part les types des
constantes de $K$ et d'autre par les types des \'{e}l\'{e}ments constants
transcendants sur $K$. Comme l'ensemble des 1-types $S_{1}(K)$ est compacte
(Proposition 3) $\langle x^{\prime}=0\rangle$ est aussi un compacte puisqu'il
est ferm\'{e} dans $S_{1}(K).$ Comme les types des constantes sont isol\'{e}s
(puisqu'ils sont r\'{e}alis\'{e}s), le type de la constante transcendante est
alors non isol\'{e} dans $S_{1}(K)$ car sinon $\langle x^{\prime}=0\rangle$
serait fini puisqu'il est compacte, ce qui contredit le fait que $C$ est
infini. Seulement la constante transcendante est dans $M(K)$, son type est
alors isol\'{e} dans $S_{1}(K)$; contradiction.\newline On d\'{e}duit alors
que $\langle x^{\prime}=0\rangle$ contient les constantes alg\'{e}briques sur
$K$ donc alg\'{e}brique sur $C,$ suppos\'{e} alg\'{e}briquement clos, toutes
les constantes sont alors dans $C$.\newline

\begin{remark}
On peut se contenter de plonger $K$ dans sa cl\^{o}ture diff\'{e}rentielle
$cl(K)$ (\cite{Poizat78} ch.II) \c{c}a nous suffit pour la d\'{e}monstration
de l'existence et de l'unicit\'{e} des extensions de Picard-Vessiot.
\end{remark}

\subsection{Unicit\'{e} des extensions de Picard-Vessiot}

\begin{lemma}
Si $\bar{u}=(u_{1},...,u_{n})$ est un syst\`{e}me fondamental de solutions,
alors $K(\bar{u})/K$ est une extension de Picard-Vessiot si et seulement si le
type de $\bar{u}$ sur $K$ est isol\'{e}.\newline
\end{lemma}

\begin{proof}
Condition n\'{e}cessaire: $K(\bar{u})/K$ est une extension de Picard-Vessiot,
on plonge $K(\bar{u})$ dans sa cl\^{o}ture diff\'{e}rentielle $M(K(\bar{u}))$
dans laquelle il n'y a pas de nouvelles constantes; de m\^{e}me on plonge $K$
dans sa cl\^{o}ture $M(K)$ qui se plonge dans $M(K(\bar{u}).$ Montrons que
$\bar{u}$ est donc dans $M(K)$.\newline On sait qu'il existe un syst\`{e}me
fondamental de solutions $\bar{v}$ dans $M(K)$ de fa\c{c}on que $K(\bar{v})$
soit une extension de Picard-Vessiot de $K$. Comme $M(K)$ est plong\'{e} dans
$M(K(\bar{u}))$ par inclusion, alors $\bar{u}=\sigma\bar{v}$ o\`{u} $\sigma$
est une matrice de constantes, donc $\bar{u}$ appartient \`{a} $M(K)$ puisque
$K,$ $M(K)$ et $M(K(\bar{u}))$ ont le m\^{e}me corps de constantes. En
conclusion le type de $\bar{u}$ sur $K$ est isol\'{e} dans $S_{n}(K)$: (car
$\bar{u}$ est dans $M(K)).$\newline Condition suffisante: Si le type de
$\bar{u}$ sur $K$ est isol\'{e} alors $K(\bar{u})/K$ est une extension
"atomique" (voir \cite{Poizat78} chapitre IV), on peut alors plonger
$K(\bar{u})$ dans la cl\^{o}ture diff\'{e}rentielle de $K.$ D'o\`{u}, dans
$K(\bar{u})$ il n'y a pas de nouvelles constantes, puisque dans la cl\^{o}ture
diff\'{e}rentielle de $K$ il n'y a pas de nouvelles constantes. Comme $\bar
{u}$ est un syst\`{e}me fondamental de solutions, $K(\bar{u})/K$ est une
extension de Picard- Vessiot. Soit $K(\bar{u})/K$ et $K(\bar{v})/K$ deux
extensions de Picard- Vessiot on peut les plonger dans la cl\^{o}ture
diff\'{e}rentielle de $K$ puisque, d'apr\`{e}s le lemme 4, les types de
$\bar{u}$ et de $\bar{v}$ sur $K$ sont isol\'{e}s dans $S_{n}(K)$. Il existe
donc un isomorphisme de la cl\^{o}ture diff\'{e}rentielle qui transforme
$K(\bar{u})$ en $K(\bar{v})$. Les extensions de Picard-Vessiot sont donc
uniques \`{a} un $K$-isomorphisme pr\`{e}s de la cl\^{o}ture diff\'{e}rentielle.
\end{proof}

\section{Groupe de galois diff\'{e}rentiel}

\begin{definition}
Soit $M$ un corps diff\'{e}rentiel et $K$ un sous corps diff\'{e}rentiel de
$M.$ Le groupe de tous les automorphismes diff\'{e}rentiels de $M$ fixant $K$
est appell\'{e} groupe de Galois diff\'{e}rentiel de l'extension $M/K$, qu'on
note $G$.
\end{definition}

Pour tout corps diff\'{e}rentiel interm\'{e}diaire $L,$ on d\'{e}finit
$L^{\prime}$ comme le sous groupe de $G$ form\'{e} par tous les
\'{e}l\'{e}ments qui fixent $L$. Pour tout sous groupe $H$ de $G$ on
d\'{e}finit $H^{\prime}$ comme \'{e}tant l'ensemble de tous les
\'{e}l\'{e}ments de $M$ qui restent fixe par $H$, $H^{\prime}$ est alors un
corps diff\'{e}rentiel interm\'{e}diaire entre $K$ et $M$.

\begin{lemma}
Soit $K$ un corps diff\'{e}rentiel de caract\'{e}ristique z\'{e}ro, soit $u$
un \'{e}l\'{e}ment dans une extension diff\'{e}rentielle de $K$ v\'{e}rifiant
$u^{\prime}=a\in K$ o\`{u} $a$ n'est pas une d\'{e}riv\'{e}e dans $K.$\newline
Alors $u$ est transcendant sur $K,\;K(u)$ est une extension de Picard- Vessiot
de $K$ et son groupe de Galois diff\'{e}rentiel est isomorphe au groupe
additif des constantes dans $K$\newline
\end{lemma}

\begin{proof}
On suppose que $u$ n'est pas transcendant sur $K,$ i.e $u$ v\'{e}rifie une
\'{e}quation polyn\^{o}miale sur $K$. Soit $u^{n}+b_{1}u^{n-1}+...+b_{n}=0$
(*) suppos\'{e}e irr\'{e}ductible.\newline En d\'{e}rivant cette \'{e}quation
on obtient $nu^{n-1}a+b_{1}^{\prime}u^{n-1}+...=0$ (car $u^{\prime}=a)$ tous
les coefficients sont donc nuls car on a suppos\'{e} que $(\ast)$ est
irr\'{e}ductible, on a alors $na+b_{1}^{\prime}=0$ ce qui entraine
$a=-b_{1}^{\prime}/n$ donc $a=(-b_{1}/n)^{\prime}$ ce qui contredit le fait
que $a$ n'est pas une d\'{e}riv\'{e}e dans $K$ ; $u$ est alors transcendant
sur $K$.\newline Comme le type de $u$ sur $K$ est isol\'{e} (car
l'\'{e}quation minimale de $u$ est d'ordre $1$ et $u$ est transcendant sur
$K$) $K(u)$ est une extension de Picard-Vessiot de $K.$\newline On remarque
que $(1,u)$ est un syst\`{e}me fondamental de solutions, en effet $(1,u)$ est
solution de l'\'{e}quation $y^{\prime\prime}-\frac{a^{\prime}}{a}y^{\prime}=0$
et $w(1,u)=a\neq0$, $a$ est non nul par hypoth\`{e}se puisque $a$ n'est pas
une d\'{e}riv\'{e}e dans $K$.\newline Soit $\sigma$ un $K$-automorphisme
diff\'{e}rentiel de $K(u)$, $\sigma u$ v\'{e}rifie alors l'\'{e}quation
diff\'{e}rentielle lin\'{e}aire $y^{\prime\prime}-\frac{a^{\prime}}%
{a}y^{\prime}=0$ et puisque $(1,u)$ est un syst\`{e}me fondamental de
solutions $\sigma u=u+c$ o\`{u} $c$ appartient \`{a} $C$ corps des constantes.
D'o\`{u}, \`{a} chaque automorphisme diff\'{e}rentiel est associ\'{e} une
constante. De m\^{e}me si $c$ est une constante, $u$ et $u+c$ sont solutions
de l'\'{e}quaiton $x^{\prime}=a.$ Ils ont alors un m\^{e}me type sur $K$,
$K(u)$ est donc isomorphe \`{a} $K(u+c)$ par un $K$-isomorphisme
diff\'{e}rentiel qui associe $u$ \`{a} $(u+c)$, mais $K(u+c)=K(u)$ car $c\in
K$; donc \`{a} chaque constante est associ\'{e} un automorphisme
diff\'{e}rentiel de $K(u)$.\newline
\end{proof}

\begin{remark}
Si $a$ est une d\'{e}riv\'{e}e dans $K,$ c'est \`{a} dire s'il existe $b\in K$
tel que $b^{\prime}=a,$ donc si $u$ est un \'{e}l\'{e}ment v\'{e}rifiant
$u^{\prime}=a=b^{\prime},$ $u$ appartient n\'{e}cessairement \`{a} $K$, $K(u)$
est alors identique \`{a} $K$; et le groupe de Galois diff\'{e}rentiel de
$K(u)=K$ est r\'{e}duit \`{a} l'identit\'{e}.

\begin{lemma}
Soit $K$ un corps diff\'{e}rentiel et $u$ un \'{e}l\'{e}ment satisfaisant
l'\'{e}quation $y^{\prime}-ay=0,$ o\`{u} $a$ appartient \`{a} $K$, on suppose
que $K(u)$ a le m\^{e}me corps de constantes que $K$. Alors $K(u)$ est une
extension de Picard-Vessiot de $K$, et son groupe de Galois diff\'{e}rentiel
est isomorphe \`{a} un sous groupe multiplicatif du groupe des constantes non
nul de $K$.\newline
\end{lemma}
\end{remark}

\begin{proof}
[Preuve]On a ici deux cas \`{a} \'{e}tudier\newline\textbf{1er Cas.} On
suppose que pour aucun entier $n$ on n'a dans $K$ de solutions non nulles de
$y^{\prime}=nay$, alors la condition
\[
\left\{
\begin{array}
[c]{ll}%
y^{\prime} & =ay\\
y & \neq0
\end{array}
\right.
\]
isole le type de $u$ sur $K$ si $u$ est non nul et v\'{e}rifie $y^{\prime
}=ay.$ En effet, si on suppose que $u$ v\'{e}rifie une \'{e}quation
polyn\^{o}miale \`{a} coefficients dans $K,$ soit $f(X)=X^{m}+\sum_{k=0}%
^{m-1}a_{k}x^{k}.$ En d\'{e}rivant cette \'{e}quation, on obtient
\[
f^{\prime}(X)=mX^{m-1}X^{\prime}+\sum_{k=0}^{m-1}[a_{k}^{\prime}%
X+ka_{k}X^{\prime}]X^{k-1}%
\]

Comme $f(u)$ et $f^{\prime}(u)$ sont suppos\'{e}s nulles on a:
\[
f^{\prime}(u)-maf(u)=\sum_{k=0}^{m-1}[a_{k}^{\prime}+kaa_{k}-maa_{k}]u^{k}%
\]
Si on suppose que $f(X)$ est une \'{e}quation de degr\'{e} minimum
v\'{e}rifi\'{e}e par $u$, on d\'{e}duit que $a_{k}^{\prime}+kaa_{k}-maa_{k}=0$
pour $k=0,1,...,$ $m-2,$ $m-1.$ D'o\`{u} $a_{k}^{\prime}=(m-k)aa_{k}$ pour
$k=0,$ $1,...,$ $m-1.$ Comme par hypoth\`{e}se cette \'{e}quation n'a pas de
solutions dans $K$ non nulles, les coefficients $a_{k}$ sont alors nulles, et
$F(X)=X^{m}.$ Ce qui entraine $u^{m}=0$, qui est absurde, il n'existe donc
aucun polyn\^{o}me \`{a} coefficients dans $K$ qui est annul\'{e} par
$u$.\newline La condition $\left\{
\begin{array}
[c]{ll}%
y^{\prime} & =ay\\
y & \neq0
\end{array}
\right.  $ isole le type de $u$ sur $K$, d'apr\`{e}s le Lemme 4, $K(u)$ est
une extension de Picard-Vessiot de $K$.\newline Pour toute constante $c\neq0,$
$cu$ est solution de l'\'{e}quation $y^{\prime}=ay$ car $(cu)^{\prime
}=cu^{\prime}=acu$, les \'{e}l\'{e}ments $u$ et $cu$ ont le m\^{e}me type sur
$K$, $K(u)$ est donc isomorphe \`{a} $K(cu)$ par un $K$-isomorphisme
diff\'{e}rentiel qui envoie $u$ sur $cu$. De m\^{e}me si $\sigma$ est un
automorphisme diff\'{e}rentiel de $K(u),$ $\sigma u$ v\'{e}rifie
l'\'{e}quation $y^{\prime}-ay=0$ on a alors $(\sigma u/u)^{\prime}=0$ donc
$\sigma u=cu$ o\`{u} $c$ est une constante non nulle.\newline Le groupe de
Galois diff\'{e}rentiel de l'extension $K(u)$ est alors isomorphe au groupe
multiplicatif des constantes non nulles de $K.$\newline\textbf{2-i\`{e}me
Cas.}

Soit $n$ le plus petit entier tel qu'il y'ait une solution non nulle dans $K$
de $y^{\prime}=nay$, soit $\beta$. On remarque que $u^{n}$ est aussi une
solution de $y^{\prime}=nay$ en effet: $(u^{n})^{\prime}=nu^{n-1}u^{\prime
}=nu^{n}a,$ car $u^{\prime}=au.$ On a donc $u^{n}=c\beta,$ o\`{u} $c$ est une
constante non nulle, $u$ est alors alg\'{e}brique sur $K$ et $X^{n}-c\beta$
est le polyn\^{o}me de plus petit degr\'{e} annul\'{e} par $u.$ En effet, s'il
existe un polyn\^{o}me, \`{a} coefficients dans $K,$ de degr\'{e} strictement
inf\'{e}rieur \`{a} $n$ annul\'{e} par $u$, soit $X^{m}+\sum_{0}^{m-1}%
a_{k}x^{k}$ avec $m<r$ on aura donc, en faisant le m\^{e}me calcul que dans le
1er cas $a_{k}^{\prime}=(m-k)aa_{k}$ pour $k=0,$ $1,,...,$ $m-1.$ Seulement
$n$ est le plus petit entier tel qu'il y ait une solution non nulle dans $K$
de $y^{\prime}=nay$, comme $m-k<n$ car $m<n,$ on d\'{e}duit que les
coefficients $a_{k}$ sont tous nuls et par ailleurs $X^{n}-c\beta$ est le
polyn\^{o}me de plus petit degr\'{e} \`{a} coefficients dans $K$ annul\'{e}
par $u$.\newline L'extension $K(u)$ est une extension de Picard-Vessiot de $K$
puisque $u$ est diff\'{e}rent de z\'{e}ro et $K(u)$ a le m\^{e}me corps de
constantes que $K$. De m\^{e}me le groupe de Galois diff\'{e}rentiel de
$K(u)/K$ est isomorphe au groupe des racines $n$-i\`{e}me de l'unit\'{e} qui
n'est qu'un sous groupe du groupe multiplicatif des constantes non nulles de
$K$. En effet, si $\sigma$ est un automorphisme diff\'{e}rentiel de $K(u)$
fixant $K$, $\sigma u$ doit \^{e}tre solution de l'\'{e}quation $X^{n}%
-c\beta=0$ donc $\sigma u$ est \'{e}gal au produit de $u$ par une racine
n'i\`{e}me de l'unit\'{e}. On remarque ici qu'une racine $n$'i\`{e}me de
l'unit\'{e} est une constante car si $\varepsilon^{n}=1,$ on a en d\'{e}rivant
$n\varepsilon^{\prime}\varepsilon^{n-1}=0,$ $\varepsilon^{\prime}$ est alors
nul.\newline De m\^{e}me pour toute racine $n$-i\`{e}me de l'unit\'{e}
$\varepsilon$, $u$ et $\varepsilon\,u$ v\'{e}rifient les m\^{e}mes
\'{e}quations. Ils ont alors un m\^{e}me type sur $K$. Il existe donc un
$K$-automorphisme diff\'{e}rentiel de $K(u)$ qui \`{a} $u$ associe
$\varepsilon\,u$\newline
\end{proof}

\section{\textbf{Groupes alg\'{e}briques des matrices}}

Soit $F$ un corps et $V$ un espace vectoriel de dimension $n$ sur $F$, i.e
l'ensemble des $n$-uples d'\'{e}l\'{e}ments de $F$. Soit $F[x_{1},..,x_{n}]$
l'anneau des polyn\^{o}mes en $n$ ind\'{e}termin\'{e}s \`{a} coefficients dans
$F.$ On d\'{e}finit alors une topologie sur $V,$ appel\'{e}e topologie de
Zariski, de la mani\`{e}re suivante, comme base de ferm\'{e}s, on prend
l'ensemble des z\'{e}ros des id\'{e}aux de $F[x_{1},\cdots,x_{n}]$, c'est un
ensemble clos par intersection quelconque et par r\'{e}union fini. Dans toute
la suite quand on parle de continuit\'{e} ou plus g\'{e}n\'{e}ralement de
topologie sur un groupe $G$, on consid\`{e}re toujours la topologie de Zariski
sur $G$ suppos\'{e} comme espace vectoriel.

\begin{lemma}
Soit $G$ un groupe muni de la topologie de Zariski. Si les multiplications
\`{a} droite et \`{a} gauche ainsi que l'application inversion sont des
hom\'{e}omorphismes de $G$ dans $G$ alors la composante connexe de
l'identit\'{e} est un sous groupe normal ferm\'{e}.

\begin{proof}
[Preuve]Si $G_{1}$ est la composante connexe de l'identit\'{e} elle est alors
ferm\'{e}e (car $\overline{G}_{1}=G_{1})$; $G_{1}^{-1}$ est aussi connexe
puisque c'est l'image continue d'un connexe, $G_{1}^{-1}$ contient $1$ alors
$G_{1}^{-1}$ est contenue dans $G_{1}.\mathcal{\ }$Si $c\in G_{1}$ alors
$cG_{1}$ est connexe et contient $c,$ d'o\`{u} $cG_{1}$ est contenu dans
$G_{1}$, ainsi on a bien $G_{1}$ est un sous groupe de $G$. Comme par
hypoth\`{e}se pour tout $x$ dans $G$, $x^{-1}G_{1}x$ est connexe (car image
par une application continue d'un connexe) et contient l'identit\'{e},
$x^{-1}G_{1}x$ est alors contenue dans $G_{1}$, ainsi $G_{1}$ est un sous
groupe normal de $G$.
\end{proof}
\end{lemma}

\begin{notation}
$K$ est un corps contenant une infinit\'{e} d'\'{e}l\'{e}ments.\newline$V$ est
un espace vectoriel de dimension fini sur $K.$\newline$\mathcal{E}$ est un
espace vectoriel des endomorphismes de $V.$

\begin{definition}
Soit $G$ un sous groupe du groupe des automorphismes de $V$, s'il existe un
ensemble $S$ de fonctions polyn\^{o}mes sur $\mathcal{E}$ tel que $G$ soit
compos\'{e} de tous les automorphismes $s$ de $V$ v\'{e}rifiant $P(s)=0$ pour
tout $P$ dans $S,$ on dit que $G$ est un groupe alg\'{e}brique
d'automorphisme, $S$ est dit ensemble de d\'{e}finition de $G$.

\begin{theorem}
Le groupe de Galois diff\'{e}rentiel d'une extension de Picard-Vessiot est un
groupe alg\'{e}brique de matrices sur le corps des constantes.

\begin{lemma}
Soit $K$ un corps diff\'{e}rentiel, $C$ corps des constantes de $K,$
$M=K(\bar{u})$ une extension de Picard-Vessiot de $K.$ Il existe un ensemble
$S$ de polyn\^{o}mes \`{a} coefficients dans $C$ tel que :\newline1)A chaque
isomorphisme diff\'{e}rentiel de $M/K$ est associ\'{e}e une matrice de
constantes satisfaisant $S.$\newline2)A chaque matrice non singuli\`{e}re
$(k_{ij})$ de constantes de $M$ satisfaisant $S$ est associ\'{e} un
isomorphisme diff\'{e}rentiel de $M/k$ envoyant $u_{i}$ sur $\sum_{j}%
k_{ij}u_{j}$.
\end{lemma}
\end{theorem}
\end{definition}

\begin{proof}
[Preuve]Soit $I$ le type de $\bar{u}$ sur $K$ (id\'{e}al diff\'{e}rentiel
premier de $K[\overline{X}]$ des \'{e}quations annul\'{e}es par $\bar{u})$.
Soit $x_{ij}$ $i,j=1,$ $2,$ $...,$ $n,$ $n^{2}$ ind\'{e}termin\'{e}es sur $M.$
L'application $y_{i}\rightarrow\sum_{j}x_{ij}u_{j}$ d\'{e}finit un
homomorphisme diff\'{e}rentiel de $K[\bar{y}]$ dans $M[x_{ij}]$, l'image de
$I$ est un id\'{e}al diff\'{e}rentiel de $M[x_{ij}]$ not\'{e} $\Delta$.

Soit $\omega_{\alpha}$ une base de l'espace vectoriel $M/C$, on \'{e}crit
chaque polyn\^{o}me de $\Delta$ dans la base $\omega_{\alpha}$, les
composantes sont des polyn\^{o}mes \`{a} coefficients constants, $S$ est
l'ensemble des composantes de tous les polyn\^{o}mes de $\Delta$.

1) Soit $\sigma$ un isomorphisme diff\'{e}rentiel de $M/k$ (dans la
cl\^{o}ture diff\'{e}rentielle de $K$).\newline$\sigma\bar{u}=(k_{ij})\bar{u}$
o\`{u} $(k_{ij})$ est une matrice de constantes, car $\bar{u}$ est un
syst\`{e}me fondamental de solutions, si on consid\`{e}re le compos\'{e} de
l'homomorphisme diff\'{e}rentiel de $K[\bar{y}]$ dans $K[\bar{u}]$ qui \`{a}
$y_{i}$ associe $u_{i}$ et de $\sigma$. L'image de $I$ par cette composition
est $0,$ car $\sigma\bar{u}$ a le m\^{e}me type que $\bar{u}$. On
consid\`{e}re maintenant l'application de $K[\bar{y}]$ dans $M[x_{ij}]$ qui
\`{a} $y_{i}$ associe $\sum_{j}x_{ij}u_{j}$, compos\'{e}e avec l'application
de $M[x_{ij}]$ dans $M$ qui \`{a} $x_{ij}$ associe $k_{ij}$. L'effet de cette
composition est le m\^{e}me que celui de l'autre composition. L'image de $I$
est alors $0$ et donc l'image de $\Delta$ par l'application de $M[x_{ij}]$
dans $M$ qui \`{a} $x_{ij}$ associe $k_{ij}$ est z\'{e}ro. Comme les $k_{ij}$
sont constants les composantes de chaque polyn\^{o}me de $\Delta$ sont nulles
pour les $k_{ij}.$
\end{proof}
\end{notation}

Pour la preuve du cas 2) on introduit le lemme suivant .

\begin{lemma}
Si $I$ et $J$ sont deux id\'{e}aux diff\'{e}rentiels premiers tels que $I$ est
contenu dans $J$ et les degr\'{e}s de transcendances de $K[\bar{a}]/K$ et
$K[\bar{b}]/K$ sont finis et \'{e}gaux. Alors $I=J$ $(\bar{a}$ et $\bar{b}$
annulent respectivement les \'{e}quations de $I$ et de $J$).

\begin{proof}
On note le degr\'{e} de transcendance de $K[\bar{a}]/K$ par $\langle K[\bar
{a}]/K\rangle.$ On a, $K[\bar{a}]=K[\bar{X}]/I$ (traduit par l'homomorphisme
diff\'{e}rentiel de $K[\bar{X}]$ dans $K[\bar{a}]$ qui \`{a} $x_{i}$ associe
$a_{i}$, le noyau de cet homomorphisme est $I$) de m\^{e}me $K[\bar{b}%
]=K[\bar{X}]/J$.

On consid\`{e}re une base de transcendance de $K[\bar{b}]/K$ ( partie
alg\'{e}briquement libre maximale), soit \{$\beta_{1},$ $...,$ $\beta_{n}\}$.
On rel\`{e}ve cette base de transcendance par l'homomorphisme de $K[\bar{a}]$
dans $K[\bar{b}]$ qui \`{a} $a_{i}$ associe $b_{i}$, en une partie
alg\'{e}briquement libre, soit \{$\alpha_{1},$ $\alpha_{2},$ $...,$
$\alpha_{n}\}$. C'est une base de transcendance de $K[\bar{a}]/K$ car les
degr\'{e}s de transcendance de $K[\bar{a}]/K$ et de $K[\bar{b}]/K$ sont finis
et \'{e}gaux \`{a} $n$. On consid\`{e}re maintenant l'homomorphisme de
$K(\bar{\alpha})[\bar{a}]$ dans $K(\bar{\beta})[\bar{b}],$ o\`{u} $\bar
{\alpha}=(\alpha_{1},\alpha_{2},...,\alpha_{n})$ et $\bar{\beta}=(\beta
_{1},\beta_{2},...,\beta_{n}),$ qui n'est autre que le compos\'{e} des
homomorphismes de $K(\bar{\alpha})[\bar{a}]$ dans $K(\bar{\beta})[\bar{a}]$
qui \`{a} $\alpha_{i}$ associe $\beta_{i}$ et de l'homomorphisme de
$K(\overline{\beta})[\bar{a}]$ dans $K(\bar{\beta})[\bar{b}]$ qui \`{a}
$a_{i}$ associe $b_{i}$. L'homomorphisme de $K(\bar{\alpha})[\bar{a}]$ dans
$K(\bar{\beta})[\bar{b}]$ est donc un isomorphisme de corps puisque c'est un
homomorphisme de corps (car $\bar{a}$ et $\bar{b}$ sont alg\'{e}briques
respectivement sur $K(\bar{\alpha})$ et $K(\bar{\beta})$ on a alors
$K(\bar{\alpha})[\bar{a}]=K(\bar{\alpha})(\bar{a})$ et $K(\bar{\beta})[\bar
{b}]=K(\bar{\beta})(\bar{b}))$ et tout homomorphisme de corps est un
isomorphisme, puisque le noyau est un id\'{e}al, il est donc soit r\'{e}duit
\`{a} l'\'{e}l\'{e}ment neutre de l'addition, soit le corps de d\'{e}part tout
entier. Ce qui veut dire que l'homomorphisme est identiquement nul, donc si on
suppose que l'homomorphisme de coprs est non identiquement nul c'est un
isomorphisme. Finalement l'isomorphisme de $K(\bar{\alpha},\bar{a})$ dans
$K(\bar{\beta},\bar{b})$ induit un isomorphisme de $K[\bar{X}]/I$ dans
$K[\bar{X}]/J$ et donc $I=J$.
\end{proof}
\end{lemma}

Dans la suite on montre le cas 2) du Lemme.

\begin{proof}
[Preuve du cas 2) du Lemme]On d\'{e}finit un homomorphisme diff\'{e}rentiel de
$K[\bar{y}]$ dans $M$ qui \`{a} $y_{i}$ associe $\sum_{j}k_{ij}u_{j}$, le
noyau $J$ de cet homomorphisme contient $I$ ( le type de $\bar{u}$ sur $K$),
car $k_{ij}$ satisfait $S$. Les conditions du lemme 9 sont satisfaites,
puisque $\bar{u}$ et $\sigma\bar{u}=(k_{ij})\bar{u}$ v\'{e}rifie une
\'{e}quation diff\'{e}rentielle lin\'{e}aire, les degr\'{e} de transcendance
de $K[\bar{u}]/K$ et de $K[\sigma\bar{u}]/K$ sont donc finis et egaux (le
d\'{e}terminant de la matrice $(k_{ij})$ est suppos\'{e} non nul), alors
$I=J$. On a donc un isomorphisme diff\'{e}rentiel de $K[\bar{u}]$ dans
$K[\sigma\bar{u}]$ qui \`{a} $u_{i}$ associe $\sum_{j}k_{ij}u_{j}$. On peut
prolonger cet isomorphisme aux corps des quotients c'est \`{a} dire \`{a}
$M=K(\bar{u})$. Ainsi, on obtient un isomorphisme diff\'{e}rentiel de $M/K$
envoyant $u_{i}$ sur $\sum_{j}k_{ij}u_{j}$.
\end{proof}

\begin{lemma}
Soit $\sigma$ un automorphisme diff\'{e}rentiel de $M/K$, si on note
$(C_{ij})$ sa matrice correspondante. Alors $W_{\sigma}=\det|C_{ij}|\times W$.

\begin{proof}
[Preuve]%
\[
\sigma u_{i}=\sum_{j}c_{ij}u_{j},\text{ }W_{\sigma}=|\sigma u_{i}%
^{(j)}|\mbox{ o\`u }\left\{
\begin{array}
[c]{ll}%
j & =0,...,n-1\\
i & =1,2,...,n
\end{array}
\right.
\]
\`{a} la j\`{e}me ligne et i\`{e}me colonne on a $\sigma u_{i}^{(j-1)}%
=\sum_{k}c_{ik}u_{k}^{(j-1)}$ la mattrice de Wronsky de $\sigma\bar{u}$ est
donc \'{e}gale au produit de la matrice de $\sigma$ par la transpos\'{e}e de
la matrice de Wronsky de $\bar{u}$, en consid\'{e}rant les d\'{e}terminants on
a bien $W_{\sigma}=|c_{ij}|W$.
\end{proof}
\end{lemma}

\begin{lemma}
Soit $M_{1}$ et $M_{2}$ deux extensions de Picard-Vessiot de $K$, si le plus
petit corps contenant $M_{1}$ et $M_{2}$ ne contient pas de nouvelles
constantes, c'est alors une extension de Picard-Vessiot de $K$.
\end{lemma}

\begin{proof}
On remarque qu'on peut plonger $M_{1}$ et $M_{2}$ dans la cl\^{o}ture
diff\'{e}rentielle de $K.$ On peut m\^{e}me plonger le corps compos\'{e}
$M_{1}\sqcup M_{2}$ dans la cl\^{o}ture diff\'{e}rentielle de $K$ (car
$M_{1}\sqcup M_{2}$ ne contient pas de nouvelles constantes). Soit $\bar
{u}_{1}$ et $\bar{u}_{2}$ les deux syst\`{e}mes fondamentaux de solutions
respectifs de $M_{1}$ et de $M_{2}$, soit $\bar{v}$ une parties de la reunion
de $\bar{u}_{1}$ et de $\bar{u}_{2}$ lin\'{e}airement ind\'{e}pendantes sur
les constantes, et tel que tous les conjugu\'{e}s de $\bar{u}_{1}$ et de
$\bar{u}_{2}$ sont combinaisons lin\'{e}aires de $\bar{v}$. On consid\`{e}re
alors l'\'{e}quation diff\'{e}rentielle lin\'{e}aire $W(X,\bar{v})/W(\bar{v})$
qui est \`{a} coefficients dans $K,$ car si $\sigma$ est un isomorphisme de la
cl\^{o}ture diff\'{e}rentielle de $K$ on a $W(\sigma X,\sigma\bar{v}%
)/W(\sigma\bar{v})=W(X,\bar{v})/W(\bar{v}),$ Lemme 10. $M_{1}\sqcup M_{2}$ est
donc une extension de Picard-Vessiot de $K$.
\end{proof}

\begin{lemma}
Une extension alg\'{e}brique normale de degr\'{e} fini est de Picard-Vessiot.
\end{lemma}

\begin{proof}
[Preuve]Soit $\bar{u}$ un g\'{e}n\'{e}rateur fini de l'extension. On
consid\`{e}re l'\'{e}quation diff\'{e}rentielle lin\'{e}aire $W(X,\bar
{u})/W(\bar{u}),$ si $\sigma$ est un isomorphisme diff\'{e}rentiel de
$M=k(\bar{u})$ c'est un automorphisme puisque l'extension est normale, et
$W(\sigma X,\sigma\bar{u})/W(\sigma\bar{u})=W(X,\bar{u})/W(\bar{u})$
l'\'{e}quation est alors \`{a} coefficients dans $K$, et $\bar{u}$ est un
syst\`{e}me fondamental de solutions de cette \'{e}quation. D'autre part $M/K$
a le m\^{e}me corps de constantes que $K$. En effet, si $a$ est une constante
dans $M$ qui n'est pas dans $K$ alors $a$ est alg\'{e}brique sur $K,$ puisque
$M/K$ est alg\'{e}brique, donc $a$ est alg\'{e}brique sur $C.$ Seulement $C$
est suppos\'{e} alg\'{e}briquement clos, $a$ est alors n\'{e}cessairement dans
le corps des constantes de $K$.
\end{proof}

\begin{theorem}
Une extension de Picard-Vessiot est normale.
\end{theorem}

\begin{proof}
[Preuve]Soit $M=K(\bar{u})$ une extension de Picard-Vessiot de $K$, $M$ peut
donc \^{e}tre plong\'{e}e dans la cl\^{o}ture diff\'{e}rentielle de $K$. Pour
tout \'{e}l\'{e}ment $z$ de $M$ qui n'est pas dans $K$ il existe donc un
isomorphisme $\sigma$ de la cl\^{o}ture diff\'{e}rentielle v\'{e}rifiant
$\sigma z\neq z$. Comme $\bar{u}=(u_{1},...,u_{n})$ est un syst\`{e}me
fondamental de solutions, $\sigma\bar{u}$ est aussi un syst\`{e}me fondamental
de solutions. on a: $\sigma u_{i}=\sum_{j}k_{ij}u_{j},$ o\`{u} $(k_{ij})$ est
une matrice de constantes, $\sigma$ est donc un automorphisme diff\'{e}rentiel
de $M/K$ et le th\'{e}or\`{e}me est ainsi d\'{e}montr\'{e}.
\end{proof}

\begin{lemma}
Si $M$ est une extension de Picard-Vessiot de $K$, $L$ corps diff\'{e}rentiel
interm\'{e}diaire entre $M$ et $K,$ alors $M$ est une extension de
Picard-Vessiot de $L$.
\end{lemma}

\begin{proof}
[Preuve]$M$ \`{a} le m\^{e}me corps de constantes que $L$ et $K,$ puisque $M$
est une extension de Picard-Vessiot de $K$. Soit $\bar{v}$ une partie de
$\bar{u}$ lin\'{e}airement ind\'{e}pendante sur $L$. On consid\`{e}re
l'\'{e}quation diff\'{e}rentielle lin\'{e}aire $W(X,\bar{v})/W(\bar{v})$. Si
$\sigma$ est un automorphisme diff\'{e}rentiel de $M/L,$ on a $W(\sigma
X,\sigma\bar{v})/W(\sigma\bar{v})=$ $W(x,\bar{v})/W(\bar{v})$. L'\'{e}quation
est donc \`{a} coefficients dans $L$.
\end{proof}

\begin{definition}
Pour tout corps diff\'{e}rentiel interm\'{e}diaire $L,$ on d\'{e}finit
$L^{\prime}$ comme \'{e}tant le sous groupe de $G$ de tous les automorphismes
diff\'{e}rentiels gardant $L$ \'{e}l\'{e}ment par \'{e}l\'{e}ment fixe.
$L^{\prime}$ est donc le groupe de Galois diff\'{e}rentiel de $M/L$. Pour tout
sous groupe $H$ de $G,$ on d\'{e}finit $H^{\prime}$ comme \'{e}tant l'ensemble
des invariants de $H$ dans $M$. $H^{\prime}$ est donc un corps
diff\'{e}rentiel interm\'{e}diaire entre $K$ et $M$. On dit qu'un groupe ou un
corps est ferm\'{e} s'il est \'{e}gal a son double prime.
\end{definition}

\begin{theorem}
Soit $M$ une extension de Picard-Vessiot de $K$, il y a une bijection entre
les corps diff\'{e}rentiels interm\'{e}diaires et les sous groupes
alg\'{e}briques du groupe de Galois diff\'{e}rentiel de $M/K$.
\end{theorem}

\begin{proof}
[Preuve]Soit $H$ un sous groupe de $G,$ on a alors $H^{\prime}=$ $InvariantH=$
$InvH,$ et $H"=G(M/InvH)$. Montrons que $H$ est Zariski dense dans
$H^{\prime\prime}.$ Pour cela, on fait une d\'{e}monstration par l'absurde. On
suppose le contraire, il existe donc un polyn\^{o}me $f$ en $n^{2}$ variable
\`{a} coefficients dans $C,$ corps des constantes, qui s'annule sur $H$ sans
toute fois s'annuler sur $H^{\prime\prime}$. On \'{e}tudie le cas d'une
\'{e}quation diff\'{e}rentielle lin\'{e}aire d'ordre $2$, la m\'{e}thode
\'{e}tant la m\^{e}me pour tout ordre.\newline Soit donc $M=K(u,v)$, la
matrice de Wronsky de $(u,v)$ \'{e}tant non singuli\`{e}re, soit $\left(
\begin{array}
[c]{ll}%
A & B\\
C & D
\end{array}
\right)  $ sa matrice inverse. Soit $y$ et $z$ deux ind\'{e}termin\'{e}es
diff\'{e}rentielles sur $M$. On pose $F(y,z)=f(Ay+By^{\prime},Az+Bz^{\prime
},Cy+Dy^{\prime},Cz+Dz^{\prime}).$ On pose $y=\sigma u$ et $z=\sigma v$ pour
$\sigma$ dans $H$. On a alors: $F(\sigma u,\sigma v)=0,$ pour tout $\sigma$
dans $H,$ car $f$ s'annule sur $H$ et non sur $H^{\prime\prime}.$ d'o\`{u}
$F(\sigma u,\sigma v)$ n'est pas nul pour tout $\sigma$ dans $H^{\prime\prime
}$. Parmis les polyn\^{o}mes de $M[y,z]$ ayant cette propri\'{e}t\'{e}, on
choisit $E$ ayant en plus le moins possible de termes. L'un des coefficients
des mon\^{o}mes de $E$ peut \^{e}tre pris \'{e}gal \`{a} 1.\newline Pour
$\tau$ dans $H$, soit $E_{\tau}$ le polyn\^{o}me obtenu en rempla\c{c}ant dans
$E$ les coefficients par leurs images par $\tau$. On a alors $E_{\tau
}(x,y)=\tau E(\tau^{-1}x,\tau^{-1}y),$ $E_{\tau}(\sigma u,\sigma v)=0$ pour
tout $\sigma$ dans $H$. De m\^{e}me $E-E_{\tau}$ \`{a} moins de monome que
$E,$ puisque l'un des coefficients de $E$ est $1$. Comme $(E-E_{\tau})(\sigma
u,\sigma v)=0,$ pour tout $\sigma$ dans $H$. Il doit \^{e}tre aussi nul pour
tout $\sigma$ dans $H^{\prime\prime}$. Si $E-E_{\tau}\not \equiv 0,$ il
existerait, donc, un \'{e}l\'{e}ment $\gamma$ de $M$ tel que $E-\sigma
(E-E_{\tau})$ ait moins de mon\^{o}mes que $E,$ et il serait nul pour tout
$\sigma$ dans $H,$ et non pour tout $\sigma$ dans $H^{\prime\prime}$.
Seulement on a suppos\'{e} que $E$ est le polyn\^{o}me le plus court ayant
cette propri\'{e}t\'{e}. On d\'{e}duit, alors, que $E-E_{\tau}\equiv0,$ ce qui
veut dire que les coefficients de $E$ sont dans $H^{\prime}$. On a, donc,
$E(\sigma u,\sigma v)=0$ pour tout $\sigma$ de $H^{\prime\prime},$ car les
coefficients de $E$ sont invariants par les \'{e}l\'{e}ments de $H^{\prime
\prime}$. $H$ est alors Zariski dense dans $H^{\prime\prime}$.\newline Soit
maintenant $L$ un corps diff\'{e}rentiel interm\'{e}diaire. Montrons que $L$
est ferm\'{e} au sens de la th\'{e}orie de Galois, c'est \`{a} dire
$L^{\prime\prime}=L$. Le lemme 13 entraine que $M/L$ est une extension de
Picard-Vessiot, $M/L$ est donc une extension normale (Th\'{e}or\`{e}me 3).
D'autre part $L$ est inclu dans $L^{\prime\prime}$, s'il existe un
\'{e}l\'{e}ment $a$ de $L^{\prime\prime}$ qui n'est pas dans $L$ on a par
d\'{e}finition de $L^{\prime\prime},$ $\sigma a=a$ pour tout automorphisme
$\sigma$ de $M/L$, comme $M/L$ est normale il existe un automorphisme $\tau$
de $M/L$ qui d\'{e}place $a$ et donc $\tau a\neq a$, il y a une contradiction
et $L^{\prime\prime}$ est alors identique \`{a} $L$.
\end{proof}

\begin{lemma}
Si $x$ est une ind\'{e}termin\'{e}e diff\'{e}rentielle, le corps des
constantes de $K(x)$ est le m\^{e}me que celui de $K$.
\end{lemma}

\begin{proof}
[Preuve]Soit $\frac{f}{g}$ un \'{e}l\'{e}ment irr\'{e}ductible constant de
$K(x)$, comme sa d\'{e}riv\'{e}e est nulle, $f^{\prime}g$ est alors identique
\`{a} $fg^{\prime}$. Soit $N$ et $M$ les ordres respectifs de $f$ et de $g$,
si on suppose que $N$ est strictement sup\'{e}rieur \`{a} $M$, le coefficient
du terme d'ordre $(N+1)$ dans l'\'{e}quation soit $g\frac{\partial f}{\partial
x^{(N)}}$ est n\'{e}cessairement nul (puisque $x$ est diff\'{e}rentiellement
transcendant sur $K$), comme $g$ est non nul $\frac{\partial f}{\partial
X^{(N)}}$ l'est obligatoirement, $f$ et $g$ sont donc dans $K$ et $\frac{f}%
{g}$ est une constante de $K$. Supposons que $N=M$ en identifiant les
mon\^{o}mes d'ordre maximal dans l'\'{e}galit\'{e} $f^{\prime}g=fg^{\prime}.$
On a, $g\frac{\partial f}{\partial x^{(N)}}=f\frac{\partial g}{\partial
x^{(N)}},$ comme $g$ n'a pas de facteurs communs avec $f,$ $g$ divise alors sa
d\'{e}riv\'{e}e par rapport \`{a} $X^{(N)},$ $\frac{\delta g}{\delta x^{(N)}}$
est donc nul et $\frac{f}{g}$ est une constante de $K$.
\end{proof}

\begin{theorem}
Il existe des extensions de Picard-Vessiot pour lesquelles le groupe de Galois
diff\'{e}rentiel est tout le groupe lin\'{e}aire $GL(n)$.
\end{theorem}

\begin{proof}
[Preuve]Soit $K_{0}$ un corps diff\'{e}rentiel, son sous-corps des constantes
$C$ est suppos\'{e} alg\'{e}briquement clos. Soit $M=K_{0}(x_{1},...,x_{n})$
le corps obtenu en adjoignant \`{a} $K_{0}$ $n$ ind\'{e}termin\'{e}es
diff\'{e}rentielles. D'apr\`{s} le lemme 14, $M$ ne contient pas de nouvelles
constantes. Soit $T$ une transformation lin\'{e}aire non singuli\`{e}re \`{a}
coefficients dans $C$, $Tx_{i}=\sum_{j}c_{ij}x_{j},$ $c_{ij}\in C$. On
d\'{e}finit $T$ sur tout $M,$ en d\'{e}crivant son action sur les
d\'{e}riv\'{e}es des ind\'{e}termin\'{e}es, $Tx_{i}^{(m)}=\sum_{j}c_{ij}%
x_{j}^{(m)}$. La transformation $T,$ ainsi d\'{e}finie, est donc un
automorphisme diff\'{e}rentielle de $M$. Si on note $K$ le sous- corps de $M$
invariant par toutes les transformations du type $T,$ $\frac{W(y,x_{1}%
,...,x_{n})}{W(x_{1},x_{2},...,x_{n})}$ est une \'{e}quation
diff\'{e}rentielle lin\'{e}aire en $y$ \`{a} coefficients dans $K$,
$(x_{1},x_{2},...,x_{n})$ est un syst\`{e}me fondamental de solution de cette
\'{e}quation, $M$ est alors une extension de Picard-Vessiot de $K$ et le
groupe de galois diff\'{e}rentiel de $M/K$ est tout le groupe lin\'{e}aire
$Gl(n)$.
\end{proof}

\paragraph{Quelques d\'{e}finitions}

Un groupe alg\'ebrique d'automorphismes est dit irr\'eductible si son id\'eal
associ\'e est premier.

On appelle dimension d'un groupe alg\'{e}brique irr\'{e}ductible
d'automorphisme le degr\'{e} de transcendance de l'anneau $\frac{C[\bar{x}%
]}{S[\bar{x}]}$ par rapport \`{a} $C$.

On appelle dimension d'un groupe alg\'ebrique quelconque la dimension du
groupe alg\'ebrique d'indice fini contenant l'\'el\'ement neutre du groupe.

Dans cette d\'{e}finition $C$ d\'{e}signe le cops des constantes.

$C[\bar{x}]$ : l'anneau des polyn\^omes en $n^{2}$ variables \`a coefficients
dans $C$.

$S[\bar{x}]$ ; l'id\'{e}al associ\'{e} au groupe alg\'{e}brique d'automorphismes.

\begin{lemma}
Soit $K$ un corps diff\'{e}rentiel, $C$ son sous-corps des constantes,
$(k_{1},k_{2},...,k_{r})$ des constantes dans un sur coprs diff\'{e}rentiel de
$K$. Si $(k_{1},k_{2},...,k_{r})$ sont alg\'{e}briquement d\'{e}pendantes sur
$K$ alors elles sont alg\'{e}briquement d\'{e}pendantes sur $C$.
\end{lemma}

\begin{proof}
[Preuve]Soit $f(k_{1},k_{2},...,k_{r})=0$ une relation polyn\^{o}miale \`{a}
coefficients dans $K.$ On consid\`{e}re une base de l'espace vectoriel $K/C,$
soit $u_{\alpha}$. En \'{e}crivant $f$ dans la base $u_{\alpha},$ on obtient
$f=\sum_{\alpha}h_{\alpha}u_{\alpha}$. Pour montrer que $h_{\alpha}%
(k_{1},k_{2},...,k_{r})=0,$ il suffit de prouver que les $u_{\alpha}$ sont
lin\'{e}airements ind\'{e}pendants dans $K(k_{1},k_{2},...,k_{r}).$ Ceci est
vrai puisque le Wronskien de toute partie fini des $u_{\alpha}$ est non nul.
Comme $h_{\alpha}$ est un polyn\^{o}me \`{a} coefficients dans $C$ alors
$(k_{1},k_{2},...,k_{r})$ sont alg\'{e}briquement d\'{e}pendantes sur $C$.
\end{proof}

\begin{lemma}
Si le degr\'{e} de transcendance de $K(\bar{u})/K$ est fini, et si $K$ est
contenu dans $L$, on peut alors \'{e}tendre le type de $\bar{u}$ sur $K$ en un
type sur $L$ et si $\bar{v}$ satisfait ce type on a en plus les degr\'{e}s de
transcendance de $K(\bar{u})/K$ et de $L(\bar{v})/L$ sont \'{e}gaux.
\end{lemma}

\begin{proof}
[Preuve]Premier cas: $\bar{u}=u_{1}$.\newline On note le degr\'{e} de
transcendance de $K(\bar{u})/K$ par $\langle K(\bar{u})/K\rangle$. On a,
$\langle K(u_{1})/K\rangle$ est fini, d'apr\`{e}s la propri\'{e}t\'{e} du
prolongment de type (voir \cite{Poizat78} ch III ). Le type de $u_{1}$ sur $K$
\`{a} un fils non deviant sur $L$, si $v_{1}$ satisfait ce type on a:
\[
\langle K(u_{1})/K\rangle=\langle L(v_{1})/L\rangle
\]
\newline Deuxi\`{e}me cas: $\bar{u}=(u_{1},u_{2})$.

Le type de $u_{1}$ \`{a} un fils non deviant sur $L$ satisfait par $v_{1}$, de
m\^{e}me le type de $u_{2}$ sur $K(u_{1})$ \`{a} un fils non deviant sur
$L(v_{1})$ satisfait par $v_{2}$, d'apr\`{e}s le th\'{e}or\`{e}me
d'aditivit\'{e} des degr\'{e}s de transcendance on a: $\langle K(u_{1}%
,u_{2})/K\rangle=\langle L(v_{1},v_{2})/L\rangle$. Si $\bar{u}=(u_{1}%
,u_{2},...,u_{n}),$ par induction, le type de $(u_{1},u_{2},...,u_{n-1})$
\`{a} un fils non deviant sur $L$ satisfait par $(v_{1},$ $v_{2},$ $....,$
$v_{n-1}),$ de m\^{e}me le type de $u_{n}$ sur $K(u_{1},u_{2},...,u_{n-1})$
\`{a} un fils non deviant sur $L(v_{1},v_{2},...,v_{n-1})$ satisfait par
$v_{n}$ le lemme est ainsi prouv\'{e}.
\end{proof}

\begin{theorem}
Soit $M$ une extension de Picard-Vessiot de $K$. La dimension du groupe de
Galois diff\'{e}rentiel de $M/K$ est \'{e}gale au degr\'{e} de transcendance
de $M/K$.
\end{theorem}

\begin{proof}
[Preuve]On remarque que le degr\'{e} de transcendance de $M=K(\bar{u})/K$ est
fini. En effet $\bar{u}=(u_{1},u_{2},...,u_{n})$ est un syst\`{e}me
fondamental de solutions d'une \'{e}quation diff\'{e}rentielle lin\'{e}aire
d'ordre $n$. Le degr\'{e} de transcendance de $K(u_{1},u_{2},...,u_{j}%
)/K(u_{1},u_{2},...,u_{j-1})$ est au plus \'{e}gal \`{a} $n$. D'apr\`{e}s le
th\'{e}or\`{e}me d'addivit\'{e} des degr\'{e}s de transcendance $\langle
K(\bar{u})/K\rangle$ est au plus \'{e}gal \`{a} $n^{2}$. Soit $L$ une
extension de $K$, si $\bar{u}$ est un syst\`{e}me fondamental de solutions
dans $L$, on peut \'{e}tendre le type de $\bar{u}$ sur $K$ en un type sur $L$,
si $\bar{v}$ satisfait ce type, on a:
\[
\langle K(\bar{u})/K\rangle=\langle L(\bar{v})/L\rangle
\]
\newline(Lemme 16)

On a alors $\bar{v}=A\bar{u}$ o\`{u} $A$ est une matrice de constantes dans
une certaine extension de $L$. Le degr\'{e} de transcendance de $L(\bar{v})$
par rapport \`{a} $L$ est donc \'{e}gal au degr\'{e} de transcendance de
$L(A)$ par rapport \`{a} $L$ qui est \`{e}gal au degr\'{e} de transcendance de
$L(A)$ par rapport \`{a} $C$ (d'apr\`{e}s le lemme 15) qui est aussi \'{e}gal
au degr\'{e} de transcendance de l'anneau $C[A\bar{X}]/S[A\bar{X}]$ qui n'est
autre que la dimension du groupe de Galois diff\'{e}rentiel de $M/K$.
\end{proof}

\begin{definition}
Soit $M$ une extension de Picard-Vessiot de $K$. On appelle cl\^{o}ture
alg\'{e}brique relative de $K$ dans $M$ la plus grande extension
alg\'{e}brique de $K$ dans $M$.
\end{definition}

\begin{corollary}
La composante connexe de l'identit\'{e} dans $G$ lui correspond la cl\^{o}ture
alg\'{e}brique relative de $K$ dans $M$.
\end{corollary}

\begin{proof}
[Preuve]La composante connxe $G_{1}$ de l'identit\'{e} dans $G$ est le plus
petit sous groupe alg\'{e}brique d'indice fini contenant l'identit\'{e} dans
$G$. La dimension de $G_{1}$ est \'{e}gale \`{a} la dimension de $G$ (voi
\cite{Chev} ch II Prop 6.4), comme la dimension de $G_{1}$ est \'{e}gale au
degr\'{e} de transcendance de l'extension $M/L$ ($M$ est une extension de
Picard-Vessiot de $L$), l'extension $L/K$ est alors alg\'{e}brique et c'est la
plus grande extension alg\'{e}brique de $K$ dans $M,$ puisque $G_{1}$ est le
plus petit sous groupe alg\'{e}brique d'indice fini contenant l'identit\'{e}
dans $G$. Dans cette d\'{e}monstration $L$ est le corps diff\'{e}rentiel
interm\'{e}diaire entre $M$ et $K$ correspondant \`{a} la composante de l'identit\'{e}.
\end{proof}

\end{document}